\documentclass[oneside]{amsart}
\usepackage[T1]{fontenc}
\usepackage[utf8x]{inputenc}
\usepackage[a4paper]{geometry}
\usepackage{float}
\usepackage{amstext}
\usepackage{amsthm}
\usepackage{amssymb}
\usepackage{graphicx}
\usepackage[unicode=true,
 bookmarks=true,bookmarksnumbered=false,bookmarksopen=false,
 breaklinks=false,pdfborder={0 0 1},backref=false,colorlinks=false]
 {hyperref}
\hypersetup{pdftitle={Исследование экстремальных свойств графов G(n,r,s)}}

\makeatletter
\numberwithin{equation}{section}
\numberwithin{figure}{section}
\theoremstyle{plain}
\newtheorem{thm}{\protect\theoremname}
\theoremstyle{plain}
\newtheorem*{lem*}{\protect\lemmaname}

\usepackage{cmap}       
\usepackage[T2C]{fontenc}
\usepackage[utf8x]{inputenc}
\usepackage[english,russian]{babel}

\usepackage{fancyhdr}
\fancypagestyle{plain}{}
\usepackage{cite}

\usepackage[figurewithin=none]{caption}

\makeatother

\providecommand{\lemmaname}{Lemma}
\providecommand{\theoremname}{Theorem}

\begin{document}
\title{$\text{New Turan-type bounds for Johnson graphs}$ }
\maketitle
\begin{center}
Dubinin Nikita Andreevich \footnote{Moscow institute of physics and technology, Dolgoprudny, Institutskiy
lane, h. 9}
\par\end{center}

\begin{center}
A new estimate is obtained for the number of edges in induced subgraphs
of Johnson graphs.
\par\end{center}

\section{Introduction}

In this paper, we consider the graph $G(n,r,s)$, whose vertices are
$r$--element subsets of the set $\{1,2,\dots,n\}$, and an edge
between two vertices is drawn if the size of the intersection of the
corresponding subsets is $s$. Another definition of the graph $G(n,r,s)$
: the graph vertices are the unit-cube vertices in $n$-dimensional
space, which have exactly $r$ ones in coordinate notation, and an
edge between two vertices is drawn when the distance between them
is $\sqrt{2(r-s)}$. It is clear that these two formulations are equivalent.
Graphs $G(n,r,s)$ are called Johnson graphs. They play a huge role
in combinatorial geometry problems (see, e.g., \cite{raigorodskii2014cliques,raigorodskii2013coloring,balogh2013coloring,pach2011combinatorial,szekely2002erdHos,larman1972realization,soifer2008mathematical,klee1991old,boltyanski2012excursions,raigorodskii2008three,raigor2007vokrug,raigorodskii2020new,sagdeev2019odnoi,ipatov2020modularnost,bobu2020odnom}),
in the coding theory (see, e.g., \cite{mac1979theory,bassalygo2000codes,raigorodskii2016combinatorial}),
in the Ramsey theory (see, e.g., \cite{nagy1972certain,graham1990ramsey,kupavskii2020theory})
etc.

In this paper, we study the extremal properties of the graph $G(n,r,s)$.
Namely, we investigate the number of edges in an arbitrary subgraph
of this graph. Notice that \emph{independent vertex set }of the $G$
graph --- such vertex subset, that no two vertices from this subset
are connected by an edge. \emph{Independence number }$\alpha(G)$
is the largest cardinality of an independent set of vertices of the
graph.

We denote by $r(W)$ the number of edges of the graph $G=(V,E)$ on
the set $W\subseteq V$. In other words, 
\[
r(W)=\left|\left\{ (x,y)\in E\mid x\in W,y\in W\right\} \right|.
\]

We also define 
\[
r(l)=\min_{\left|W\right|=l,W\subseteq V}r(W).
\]

The question arises about the study of this value. Classical Turan's
theorem 1941 gives the answer to this question in the general case.
\begin{thm}
Let $G$ be an arbitrary graph, let $\alpha$ be its independence
number, $l>\alpha$ . Then $r(l)\geqslant\frac{l^{2}}{2a}-\frac{l}{2}$.
\end{thm}

The proof of this theorem does not take into account any special properties
of the graph $G$, and, moreover, this theorem is not improveable
in general. However, it is reasonable to assume that for graphs with
some constraints, the estimate can be improved. We consider \emph{distant}
graphs --- graphs, whose vertices are points in $\mathbb{R}^{n}$
space, and the edge between such vertices is present if and only if
the distance between them is equal to some constant. It is clear that
defined graph $G(n,r,s)$ is distant graph.

For arbitrary distant graphs the following theorem has been proven(\negthinspace{}\negthinspace{}\cite{mihailov2009chislah}).
\begin{thm}
Let $G_{n}$ be sequence of the distant graphs, which $V(G_{n})\subset\mathbb{R}^{n}$.
Let $\alpha_{n}=\alpha(G_{n})$. Let $W_{n}$ be a subset of $V(G_{n})$.
If$\left|W_{n}\right|=l(n)$ and $n\alpha_{n}=o(l(n))$, then with
$n\rightarrow\infty$

\[
r(l(n))\geqslant\frac{l(n)^{2}}{\alpha_{n}}(1+o(1)).
\]

\medskip{}
\end{thm}

Then, we see, that on distant graphs, forming sequences with certain
asymptotic properties, Turan assessment has been improved twice. It
can be assumed that at an even more narrow class of distant graphs
$G(n,r,s)$ the assessment allows for further improvements. And indeed,
in the \cite{pushnyakovchisle} paper the following theorem has been
proven (see, e.g., \cite{pushnyak2015new,pushnyak2015chisle,pushnyakov2018otsenka,pushnyakovkolichestvah}).
\begin{thm}
Consider graph $G(n,3,1).$ Let the function $l:\mathbb{N}\mathbb{\rightarrow N}$
satisfy $n^{2}=o(l)$ as $n\rightarrow\infty$. Then there is a function
$h:\mathbb{N}\mathbb{\rightarrow N}$, that $h\sim\frac{3l^{2}}{2n}$
with $n\rightarrow\infty$ and $r(l(n))\geqslant h(n)$ for any large
enough $n\in\mathbb{N}.$
\end{thm}

To understand how the results of theorems 2 and 3 correlate, note,
that $\ensuremath{\alpha(G(n,3,1))\in\{n-2,n-1,n\}}$ (see \cite{pushnyakovchisle,nagy1972certain}).
This means, that on its class of graphs theorem 3 is one and a half
times stronger than the general theorem 2.Our main result will be
a generalization of theorem 3 in case of fixed $r,s$ with the condition,
that $r=2s+1$ and $r-s$ is the power of prime number. Obviously,
the parameters of theorem 3 satisfy these conditions. Note, that it
is in this conditions in the \cite{frankl1985forbidding} was shown,
that $\alpha(G(n,r,s))\sim n^{s}\frac{(2r-2s-1)!}{r!(r-s-1)!}$. This
means that there is a function $q(n)=(1+o(1))$, that $\alpha(G(n,r,s))=q(n)\cdot n^{s}\frac{(2r-2s-1)!}{r!(r-s-1)!}.$
Also, since the function $q(n)$ is limited, there is the constant
$C_{0},\text{ that }\alpha(G(n,r,s))\leqslant C_{0}n^{s}$. We will
sometimes need these formulations. So, we have
\begin{thm}
Let $r=2s+1$ and $r-s$ is power of the prime number, and $l(n)$
is any function with limitations $l(n)=o(n^{2s+1})$ and $n^{2s}=o(l(n))$.
Let $\alpha_{n}=\alpha(G(n,r,s))$. Then there is function $h:\mathbb{N}\rightarrow\mathbb{N},$
satisfy  $h\sim\frac{3l(n)^{2}}{2\alpha_{n}}$ with $n\rightarrow\infty$,
and $r(l(n))\geqslant h(n)$.
\end{thm}

To prove theorem 4, we will need an additional lemma.
\begin{lem*}
Let the parameters $r,s$ and function $l$ satisfy the conditions
of the theorem 4. Let $W$ be a random set of the vertices of the
graph $G(n,r,s)$, whose size is $l(n)$. Let $\Gamma$ be a maximal
independent set of vertices in a subgraph of a graph $G(n,r,s)$,
based on vertices from $W$. Let $w\in W\setminus\Gamma$. Denote
by $n(\Gamma,w)$ number of vertices in $\Gamma$, related to the
$w$. Let $U_{1}$ and $U_{2}$ be  sets of such vertices $w\in W\setminus\Gamma,$
that $n(\Gamma\mathit{,w\mathrm{)}}=1$ or 2 respectively. Then there
exists a constant $C_{1}$ such that $\left|U_{1}\cup U_{2}\right|\leqslant C_{1}n^{2s}$.
\end{lem*}
We emphasize that the constant in the lemma will depend only on $r$
and $s$, but not on the $l$, $W$ or $\Gamma$.

In the next section, we first give a proof of the lemma, and then,
in subsection 2.2, we prove theorem 4. In the proofs, in order to
avoid confusion, it will sometimes be convenient for us to distinguish
between the notation for one vertex $u$ or another of the graph $G(n,r,s)$
and the corresponding $r$-element subset. The latter will be denoted
by $\text{supp}(u)$ and will be called the \emph{support} of the
vertex $u$.

Finally, note that similar results for the case of arbitrary distance
graphs in the plane can be found in the paper \cite{shabanov2016turan}

\section{Proofs}

\subsection{Proof of the lemma}

Firstly, let us prove that there exists $C_{2}$such that $\left|U_{1}\right|\leqslant C_{2}n^{s+1}$.
Let's choose the vertex $u\in\Gamma.$ Let
\[
U_{1,u}=\left\{ w:w\in W\setminus\Gamma,n(\Gamma,w)=1\text{ and }(w,u)\in E\right\} .
\]
\begin{figure}[H]
\caption{Vertices with single edge with $\Gamma$}

\includegraphics[width=5cm,height=5cm]{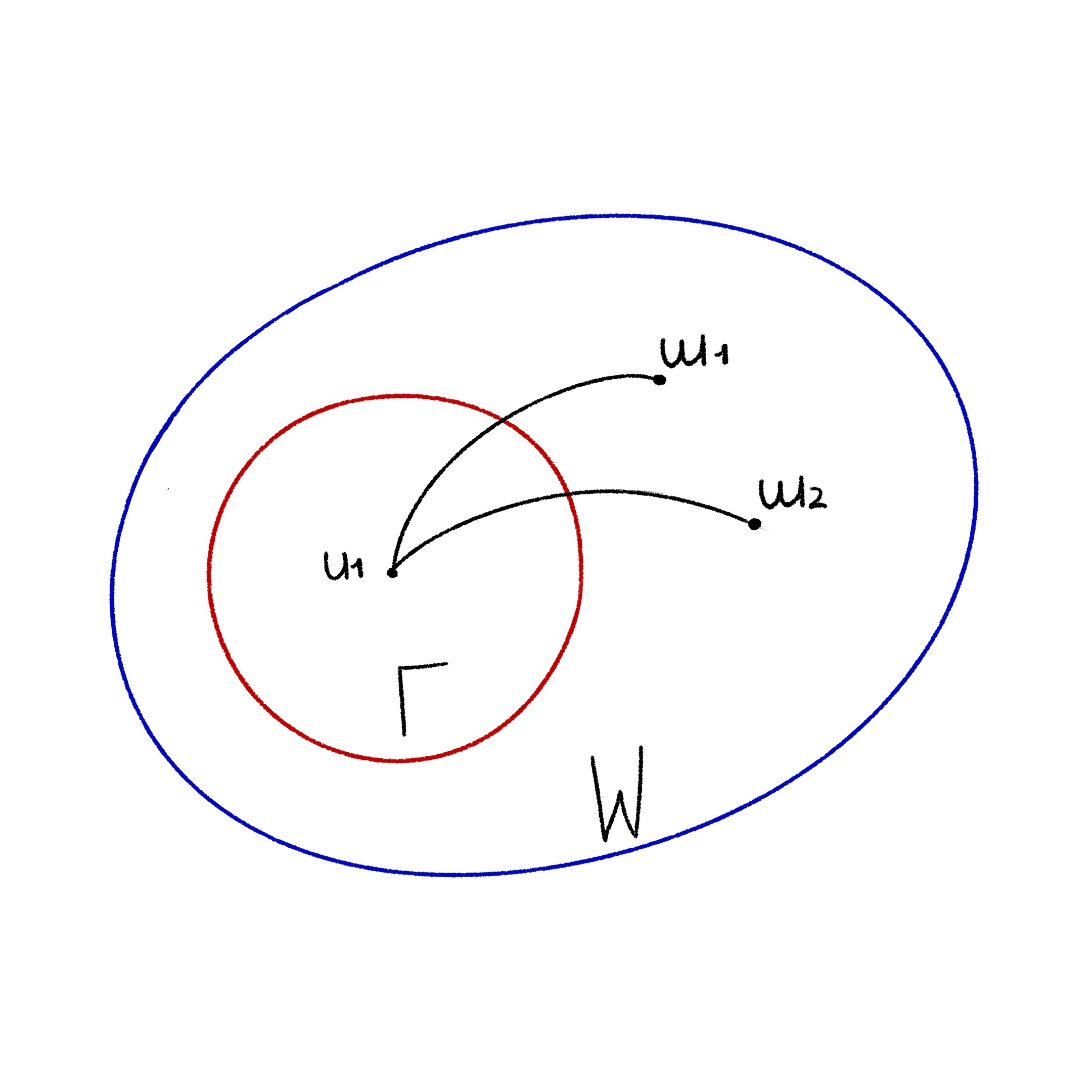}
\end{figure}
Denote by $U_{1}$ the union of sets $U_{1,u}$ by all $u\in\Gamma$.
Let us fix $u$ and estimate the cardinality of $U_{1,u}.$ Let $v\in U_{1,u}$.
supports of the vertices $u\text{ and }v$ intersect by $s$ elements,
and these elements can be selected in $C_{r}^{s}$ ways. Next, $(s+1)$,
element of the support $v$ can be chosen in $n-r$ ways. And for
its $s+1$ elements there is at most one way to select all the others.
Let it is wrong, so there are at least two vertices $v_{1}=\left\{ v^{1},v^{2},\dots,v^{s+1},a,\dots\right\} $
and $v_{2}=\left\{ v^{1},v^{2},\dots,v^{s+1},b,\dots\right\} $. Note
that there is no edge between them, because its supports intersect
by $s+1$ elements. Also note, that each of the vertices has only
one edge with the set $\Gamma$ and this edge leads to the selected
vertex $u$ (see pic. 1). Then the set $(\Gamma\setminus\left\{ u\right\} )\cup\left\{ v_{1},v_{2}\right\} $
doesn't have any edges, that is, it is independent, and has cardinality
greater than the cardinality of the set $\Gamma$, which contradicts
the assumption of maximality of $\Gamma$. So, 
\[
\qquad\left|U_{1}\right|\leq C_{r}^{s}\cdot\left|\Gamma\right|\cdot n\leqslant C_{r}^{s}\cdot\alpha_{n}\cdot n\sim n^{s+1}C_{r}^{s}\frac{(2r-2s-1)!}{r!(r-s-1)!}.
\]

Now let us prove that there is such a constant $C_{3}\text{, that }\left|U_{2}\right|\leqslant C_{3}n^{2s}$.
In other words, it is necessary to estimate the number of such vertices
$w\in W\setminus\Gamma,$ that $n(\Gamma\mathit{,w\mathrm{)}}=2.$
Let $w$ have the edges with $u_{1},u_{2}\in\Gamma.$ supports of
the $u_{1},u_{2}$ can intersect in $0,1,\ldots,s-1,s+1,\ldots,r-1=2s$
elements. Since the vertex $w$ has an edge with any vertex $u_{1},u_{2}$,
support of $w$ and union of supports of $u_{1},u_{2}$ can intersect
in $s,s+1,\ldots,r-1=2s$ elements. Let us define\emph{ checkmark.
}We will call \emph{checkmark }three vertices with such properties:
$u_{1},u_{2}\in\Gamma$ and $w\in W\setminus\Gamma,$ $n(\Gamma\mathit{,w\mathrm{)}}=2$,
moreover $(u_{1},w)\in E,(u_{2},w)\in E$ (see pic. 2). Vertex $w$
will be called the\emph{ center }of the checkmark, other checkmark
vertices will be called \emph{sides.} We divide the proof into two
cases. In the first case support of the center of the checkmark intersects
with union of the sides' supports more, than by $s$ elements. In
the second case checkmark's center's support intersects with union
sides' supports by $s$ elements. It is possible if supports of the
checkmark's sides have $s+1$ common element (not $s$ elements, because
$u_{1},u_{2}$ don't have an edge).
\begin{figure}[H]
\includegraphics[width=5cm]{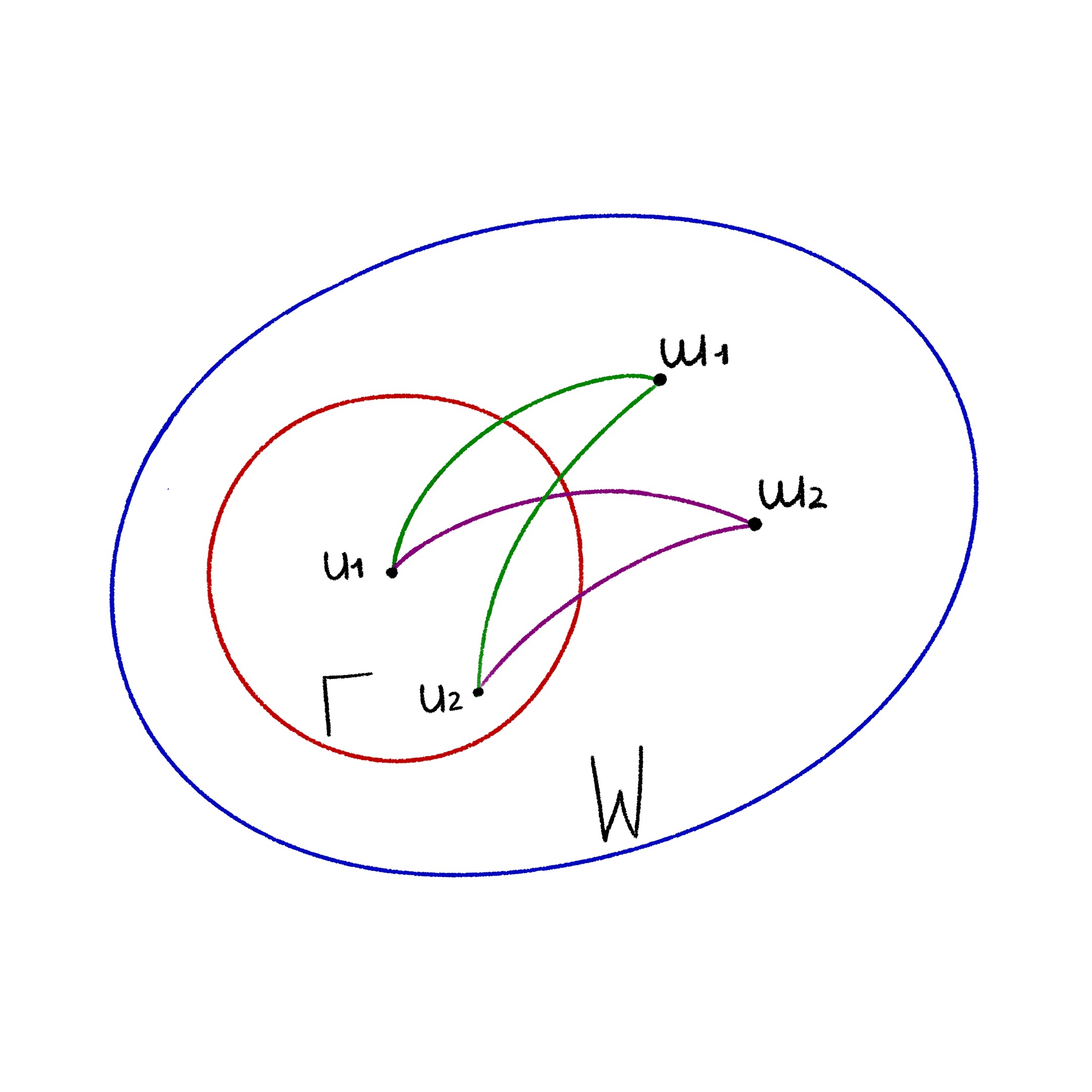}

\caption{Checkmarks}
\end{figure}

\begin{figure}[H]
\includegraphics[width=5cm]{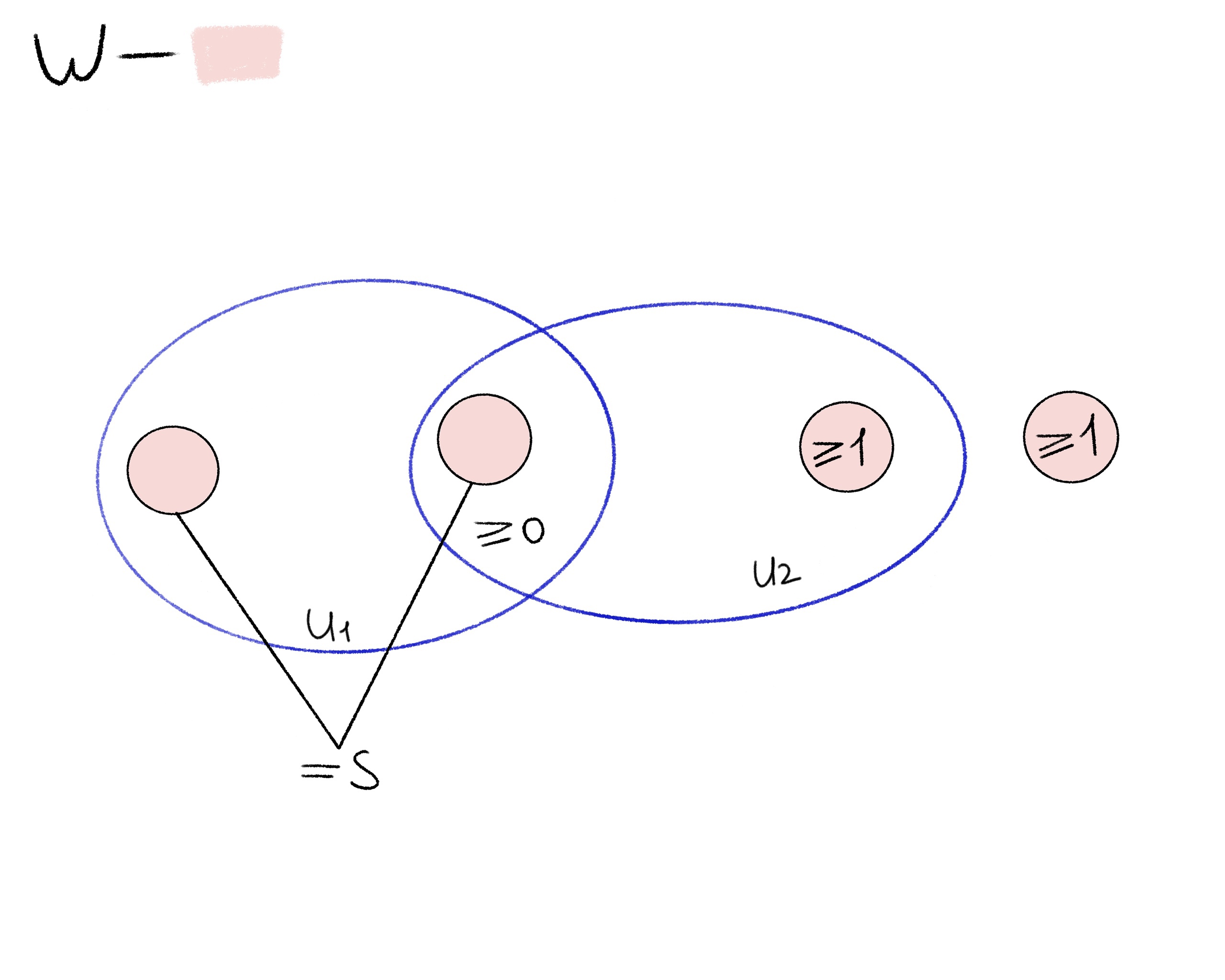}\caption{supports in the case 1}

\end{figure}
\emph{Case 1}. Let$\left|(\text{supp}(u_{1})\cup\text{supp}(u_{2}))\cap\text{supp}(w)\right|=k,k\geq s+1.$
Choose pair $u_{1},u_{2}\in\Gamma$ and count number of the checkmarks
with sides in this pair and centers in certain $w$. As long as $\left|\text{supp}(u_{1})\cap\text{supp}(w)\right|=s,$
there are $C_{r}^{s}$ ways choose $s$ elements of the support $w$,
by which it will intersect with the support $u_{1}$. Among not more
than $r$ elements in the difference of the support$u_{2}$ and support
of $u_{1}$, we have to choose $k-s\geqslant1$ elements (see pic.
3). So, $s+1$-th element is chosen in at most r ways in the support
$u_{2}$. And there are no more than two ways to select the remaining
elements, and it doesn't matter where to choose these elements ---
in the support or in the whole set of $n$ numbers. Suppose the contrary,
then there are at least three vertices $w_{1},w_{2},w_{3}$, whose
supports have $s+1$ common element, and therefore, there are no edges
between them. Therefore, since from each vertex $w_{1},w_{2},w_{3}$
there are exactly 2 edges with the set $\Gamma$, and all the edges
go to the vertices $u_{1},u_{2}$, the set $(\Gamma\setminus\{u_{1},u_{2}\})\cup\{w_{1},w_{2},w_{3}\}$
is independent and has a larger cardinality, than $\Gamma$, which
contradicts the maximality of $\Gamma$. So, for any two vertices
$u_{1},u_{2}$ there are no more $2r\cdot C_{r}^{s}$ checkmarks.
A pair of vertices can be selected at most $\alpha_{n}^{2}\sim n^{2s}\left(\frac{(2r-2s-1)!}{r!\cdot(r-s-1)!}\right)^{2}$
ways, and therefore, the number of such vertices $w$ does not exceed
\[
\qquad\qquad\qquad\alpha_{n}^{2}\cdot2r\cdot C_{r}^{s}<C_{4}n^{2s}.
\]
\begin{figure}[H]
\includegraphics[width=5cm]{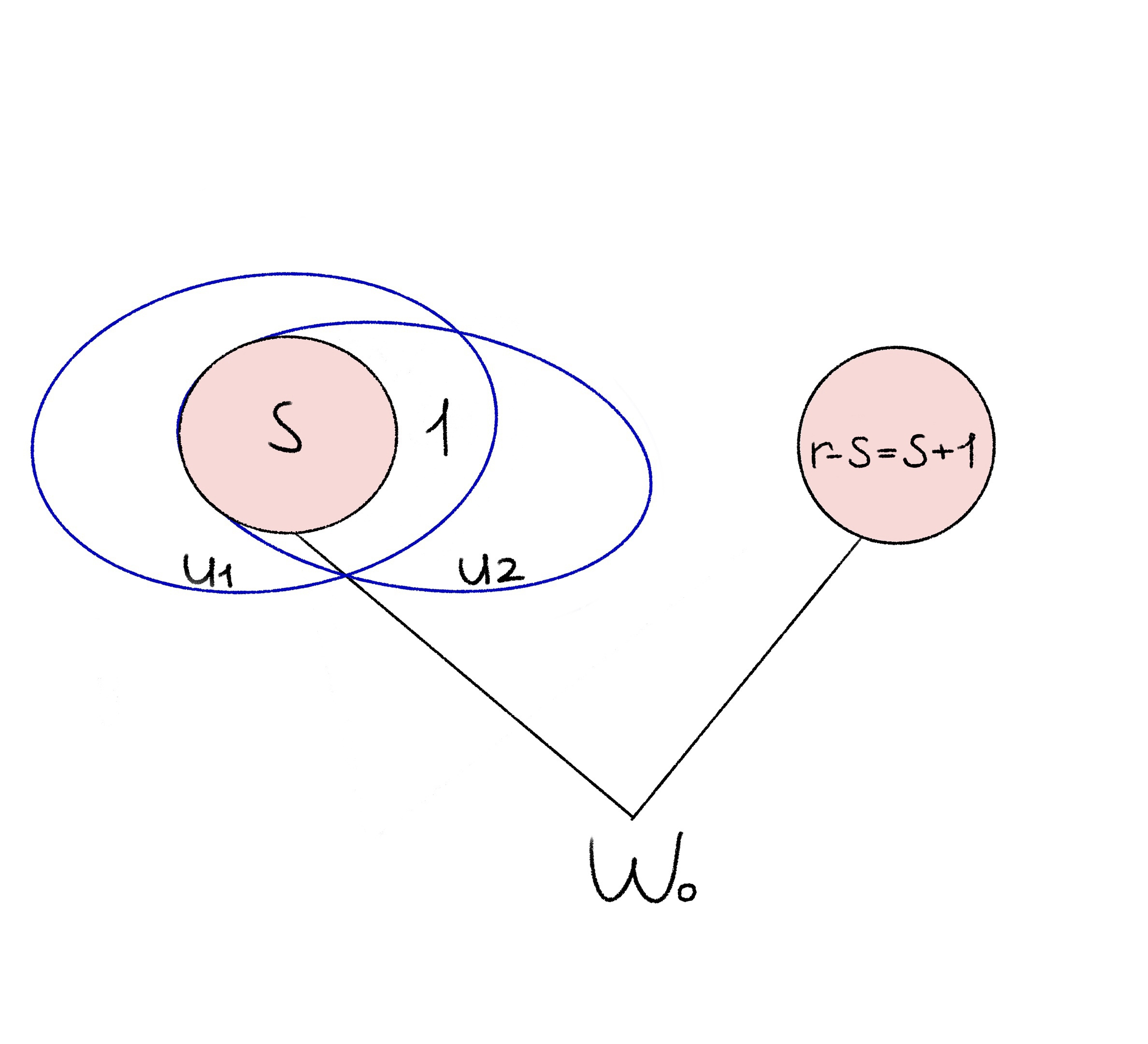}\caption{supports in the case 2}

\end{figure}
\emph{Case 2}. Here we are interested in checkmarks that have$\left|(\text{supp}(u_{1})\cup\text{supp}(u_{2}))\cap\text{supp}(w)\right|=s$.
Choose one vertex $u_{1}\in\Gamma$, fix $s$ elements in its support.
Let there be at least one such checkmark that the support of its center
$w_{0}$ intersects with its sides supports $u_{1},u_{2}$ precisely
along these $s$ elements (we will call them \emph{fixed}), and its
sides supports, respectively, intersect by fixed $s$ and at least
one more element (see pic. 4). Now let's calculate how many more sides
of the checkmarks with the side at the vertex $u_{1}$ and the same
$s$ elements (we want to estimate exactly the number of sides; we
will estimate the number of centers for a given pair of sides later).
Since all supports of the sides of the checkmarks have common fixed
elements with the support $u_{1}$ and there are $s$ such elements,
each of the side supports must have some $(s+1)$-th common element
with support of the $u_{1}$, otherwise, in an independent set of
vertices $\Gamma$ an edge is formed. Moreover, all checkmark sides,
except $u_{1}$ and the second checkmark side with the center in $w_{0}$,
do not have edges with $w_{0}$, because the center of the checkmark
has exactly two edges with the set $\Gamma$. However, all the supports
of these sides and the support $w_{0}$ have the same fixed $s$ elements
in the intersection. This means each of these sides supports must
have at least one additional common element with the support $w_{0}$.
Thus, the required number of sides is at most $(r-s)^{2}C_{n}^{r-s-2}\le C_{5}n^{r-s-2}=C_{5}n^{s-1}$.
Now, for each pair of sides of the checkmark, we will count how many
centers $w$ can exist. First $s$ elements of the support are fixed
--- $w$ intersects with supports of $u_{1}$and $u_{2}$ by this
fixed elements. Next element we choose not more that $n$ ways besides
the supports $u_{1}$ and $u_{2}$. And there are no more than two
ways to select all other elements. The proof is similar to the proof
in the case 1 --- otherwise we will have three vertices $w_{1},w_{2},w_{3}$
without an edges, and the set $(\Gamma\setminus\{u_{1},u_{2}\})\cup\{w_{1},w_{2},w_{3}\}$
will have cardinality more than $\Gamma$, and still be independent.
Hence, for the vertex $u_{1}$ and some $s$ elements from its support
checkmarks centered at some $w$ at most $2n$, whence we get that
the number of checkmarks with the given $u_{1}$ and given $s$ fixed
elements no more than $2n\cdot C_{5}n^{s-1}=C_{6}n^{s}$. Ways to
choose the $u_{1}$ vertex and $s$ elements in it respectively 
\[
\qquad\qquad\qquad\alpha_{n}\cdot C_{r}^{s}\leqslant n^{s}C_{7},
\]
hence, in the current case, the number of vertices $w$ is at most
\[
\qquad\qquad\qquad(C_{6}n^{s})\cdot(C_{7}n^{s})=C_{8}n^{2s}.
\]
So, in each of the cases, we have an estimate of the form $C\cdot n^{2s}$.
Adding all the constants, we obtain the value stated in the lemma
$C_{1}n^{2s}$ .

\subsection{Proof of theorem 4}

Let $W$ --- some subset in the set of vertices of the graph $G(n,r,s)$,
having cardinality $l=l(n)$. Consider the largest independent set
in terms of cardinality $\Gamma_{1}$ in the subgraph of $G(n,r,s)$,
generated by the set of vertices $W$. Let its cardinality equal $\beta_{1}\leqslant\alpha_{n}$.
Let $F_{1}$ --- a subset of such vertices in the set $W\setminus\Gamma_{1}$,
that for any vertex $w\in F_{1},$ $n(\Gamma_{1},w)\leqslant2.$ Let
$f_{1}=\left|F_{1}\right|$. It follows from the lemma that $f_{1}\leqslant C_{1}n^{2s}$.
Note that any vertex $u\in W\setminus(\Gamma_{1}\cup F_{1})$ has
at least three edges with $\Gamma_{1}$. Then found at least $3(l(n)-f_{1}-\beta_{1})+f_{1}\geqslant3(l(n)-\alpha_{n})-2C_{1}n^{2s}$
edges. Let us remove from $W$ the independent set$\Gamma_{1}$ and
in the resulting set $W\setminus\Gamma_{1}$ choose a new largest
independent set $\Gamma_{2}$ with cardinality equal $\beta_{2}\leqslant\alpha_{n}$.
Let $F_{2}$ --- such subset of the set $W\setminus(\Gamma_{1}\cup\Gamma_{2})$,
that for any vertex $w\in F_{2}$ we have $n(\Gamma_{2},w)\leqslant2.$
Let $f_{2}=\left|F_{2}\right|$. From the lemma we have the estimate
$f_{2}\leqslant C_{1}n^{2s}$. We found again at least $3(l(n)-2\alpha_{n})-2C_{1}n^{2s}$
edges. Repeat this operation $\left[\frac{l(n)}{a_{n}}\right]$ times,
we get an estimate 
\[
\begin{aligned}r(l(n))\geqslant\sum_{i=1}^{\left[l(n)/\alpha_{n}\right]}\left(3(l(n)-i\alpha_{n})-2C_{1}n^{2s}\right)\sim3l(n)\cdot\frac{l(n)}{\alpha_{n}}-\frac{3}{2}\alpha_{n}\cdot\frac{l(n)}{a_{n}}\left(\frac{l(n)}{a_{n}}+1\right)-\frac{l(n)}{\alpha_{n}}\cdot2C_{1}n^{2s}\sim\\
\sim3\frac{l^{2}(n)}{\alpha_{n}}-\frac{3}{2}\frac{l^{2}(n)}{a_{n}}-2C_{1}n^{2s}\frac{l(n)}{a_{n}}\sim\frac{3}{2}\frac{l^{2}(n)}{a_{n}}
\end{aligned}
\]
 with $n\rightarrow\infty$, since under the assumption $n^{2s}=o(l(n))$.
Theorem 4 is proven.

\vspace{0.5cm}

The author is grateful to Andrei Mikhailovich Raigorodsky for his
multifaceted support, without which the work would not have taken
place.

The author is also grateful and expresses gratitude to the artist
of drawings and diagrams, HSE student Maria Smetanina.

\bibliographystyle{plain}
\bibliography{bibtex}

\end{document}